\documentclass[francais]{article-arima-tamtam}
\author{Hamad Hazim\fup{*}\andauthor Bernard Rousselet\fup{**}}

\title{Frequency sweep for a beam system with}
\subtitle{local unilateral contact modeling satellite solar arrays}
\address{\fup{*}
  U.N.S.,   Laboratoire J.A. Dieudonn\'e \\
U.M.R.  C.N.R.S. 6621, Parc Valrose, \\ F 06108 Nice,
 C{\'e}dex 2, email~:
\\
hamad.hazim@unice.fr\\
\fup{**}
  U.N.S.,   Laboratoire J.A. Dieudonn\'e \\
U.M.R.  C.N.R.S. 6621, Parc Valrose, \\ F 06108 Nice,
 C{\'e}dex 2, email~: {\sf br\char64math.unice.fr}
 }%
\resume{Dans le but  de r\'eduire la masse des panneaux solaires
d'un satellite, la flexibilit\'e de ces panneaux  est devenue non
n\'egligeable et ils peuvent se taper l'un contre l'autre ce qui peut provoquer un endomagement de la structure.
Pour prevenir cela, des cales sont fix\'{e}es dans des positions bien choisies de la structure; ils
jouent le r\^{o}le d'un ressort \'{e}lastique unilateral. Comme cons\'equence n\'egative, la dynamique de ces panneaux devient
nonlin\'eaire. Une approximation par \'elements finis est utilis\'ee
pour r\'esoudre les \'equations aux d\'eriv\'ees partielles qui gouvernent la dynamique de la structure;
des balayages en fr\'equence ont \'et\'e faits
num\'eriquement pour \'etudier le comportement dynamique. Les modes normaux non lin\'aires sont \`{a} l'\'etude.}
\abstract{In order to save mass of  satellite solar arrays, the flexibility of
the panels becomes not negligible and they may strike each other;
this may damage the structure. To prevent this, rubber snubbers are
mounted at well chosen points of the structure and they act as \textbf{one sided} linear
spring; as a negative consequence, the dynamic of these panels becomes \textbf{nonlinear}.
 The finite element approximation is
used to solve partial differential equations governing the structural
dynamic. Frequency sweep has been performed numerically to study the
dynamic behavior. Non linear normal modes are under study}
\motscles{Dynamique non lin\'eaire, Syst\`emes differentiels, Mod\'elisation}
\journal[TAM-TAM, Maroc, Kenitra 6-8 may 2009]
{\textbf{TAM-TAM'09}}{1}{1}{2003}{1}{7}

\setboolean{isrubrique}{true}
\begin{document}
\maketitlepage
%
\section{Introduction}

In articles $\cite{r1}$, $\cite{r2}$, $\cite{r4}$ and $\cite{r5}$, the dynamic of a beam system with a nonlinear contact force, under a periodic excitation given as an imposed acceleration form is studied both numerically and experimentally. When sweeping frequencies in an interval which contains eigen frequencies of the beam, resonance phenomena appear as well as new frequencies caused by the unilateral contact.\\

\begin{figure}[hbtp]
\includegraphics[scale=.24]{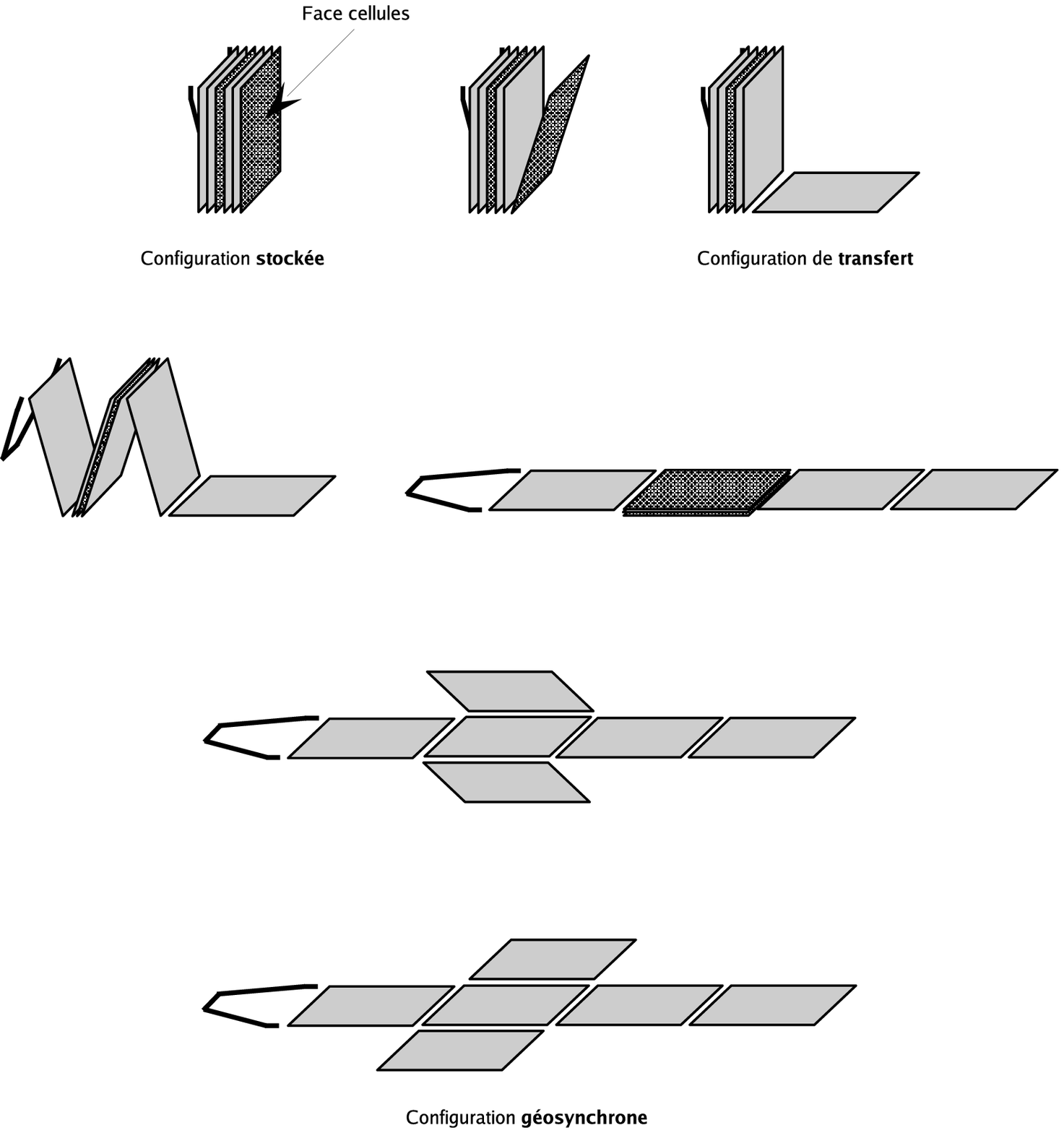}
\hspace{1.0cm}
\includegraphics[scale=.24]{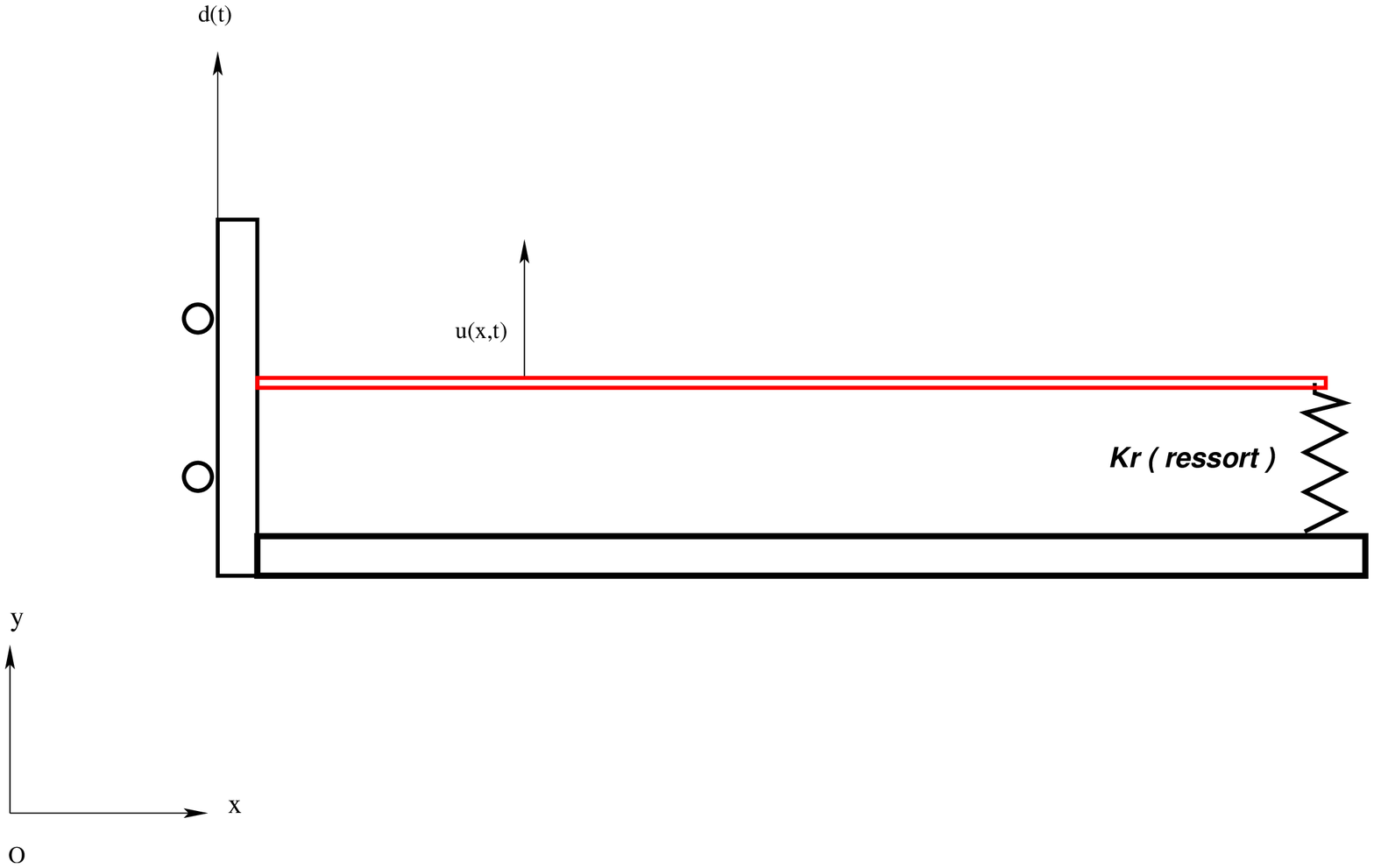}
\caption{\small{At left: solar arrays from folded to final position,
at right: beam system with unilateral spring on a shaker to model an unfolded array}} \label{1}
\end{figure}

The study of the total  dynamic behavior of solar arrays in a folded position with snubbers  are so complicated, that to simplify, a solar array is modeled by a clamped-free Bernoulli beam with one-sided linear spring. This system is fixed on a shaker which has a vibratory motion $d(t)$, it is given as a
periodic imposed acceleration: $\ddot{d}(t)=a\sin\omega t$ (see figure $\ref{1}$).\\

The present study is to simulate the behavior of a beam which strickes a snuber. In this first step, existing algorithms are used in order to prepare a sequence of experiments.\\
The focus is on the comparison between the linear and the non-linear case.\\
An original approach of an extension of normal modes to this non-linear case is under study;  this approach is an alternative to the one of
\cite{a3} and \cite{a4}\\
The beam motion with a snubber can be modeled as:

\begin{equation}
\rho
S\ddot{u}(x,t)+EIu^{(iv)}(x,t)=0
\end{equation}

with the boundary conditions:\\
$u(0,t)=d(t)$, $\partial_x u(0,t)=0,$\\
$EIu^{(2)}(L,t)=0$ and $EIu^{(3)}(L,t)=k_{r}(d(t)-u(L,t))_{+}$.\\
\begin{equation}
u(x,t)_+=\left\lbrace
\begin{array}{rl}
u(x,t) & \ if \ u > 0\\
0  & \ if \ u \leq 0\\
\end{array} \right.
\end{equation}
The classical Hermite cubic finite element approximation is used, it yields an ordinary differential system in the form:
\begin{equation}
M\ddot{q}+Kq=k_{r}(d(t)-q_{n})_+e_{n}
\end{equation}
Where $M$ and $K$ are respectively the mass and stiffness assembled matrices, q is the vector of degrees of freedom of the beam,
$q_i=(u_i,\partial_xu_i)$, $i=1,...,n$, where $n$ is the size of $M$.\\
The system is integrated numerically using the Scilab routine 'ODE' for  stiff problems, package ODEPACK is called and it uses the BDF method $\cite{a6}$.\\

The frequency sweep is performed by running the integration of ordinary
differential systems at a chosen initial value of the frequency and
saving the maximum of the displacement of all nodes in the whole
time interval. The frequency is then incremented till the end of the
frequency interval and the same procedure is used again.\\
The maximum of displacement is plotted in the frequency domain, the simulation
shows many resonance phenomena when approaching  particular frequencies;
these particular frequencies are the non-linear frequencies of the system corresponding to
the eigen frequencies of the linearized system (bilateral spring). Other superharmonic and subharmonic
frequencies appear due to the nonlinearity of the contact.
The frequency responses are compared with the FFT of the temporal signal coming from
the direct integration of the differential systems. On the other hand, the calculation of the \textbf{non-linear normal
modes (NNM)} \cite{a4} and \cite{a3} will be a theoretical way to validate the numerical simulations.

\section{Numerical results}

The Physical properties and the dimensions of the beam with a spring stiffness of $5\times 10^5 N/m$ gives the following eigen-frequencies (bilateral spring):\\

\begin{itemize}
\item $first\; frequency= 196.35694  Hz$ \hspace{1cm} $second\;frequency= 472.07584  Hz$
\item $third \;frequency=961.52462 Hz$
\end{itemize}
The figure $(\ref{f1})$ show the eigen-frequencies of the linearised system (bilateral spring) calculated analytically and the corresponding frequencies of the non-linear cases (unilateral spring);  many other frequencies appear in the non-linear cases due to the nonlinearity of the contact. The superharmonics and the subharmonics frequencies
are shown too, the superharmonics are the double and the triple (etc..) of the corresponding eigen frequency, the subharmonics are the half and the third (etc..) of the eigen frequencies. Many other combinations of these frequencies may exist but it is quit difficult to distinguish their location.\\
A damping term in the spring will help to distinguish each mode with a particular initial condition or excitation, the full  article will show this technique.\\
The numerical results are to be compared with experiments in preparation.\\
On the other hand, the calculation of the non-linear normal modes are under development in order to address more precisely the frequency response of the structure and to validate numerical results.\\
Comparison with results obtained by asymptotic expansions of the type in $\cite{JB}$ is in project
\begin{figure}[hbtp]
 \includegraphics[width=16cm]{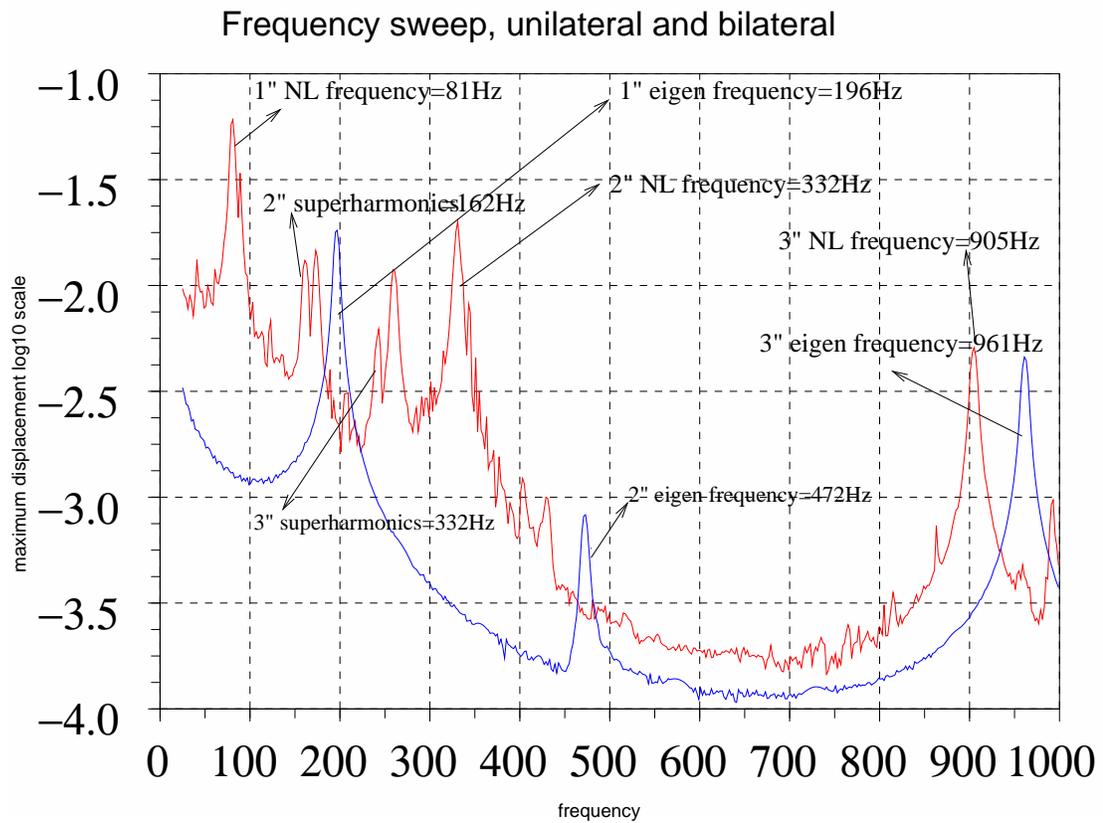}
\caption{Frequency sweep, bilateral and unilateral spring, no pres stress or backlash, $ beam \ height=10mm$, $kr=500.000 N/m$, $a=50m/s^2$, $tf=0.4s$}
\end{figure}\label{f1}

\newpage

\end{document}